\documentclass[12pt]{article}
\usepackage{amsfonts}
\usepackage{amsmath}
\usepackage{verbatim}
\usepackage{amssymb}
\usepackage[latin1]{inputenc}
\usepackage{a4wide}
\usepackage{amsthm}
\usepackage{graphicx}

\DeclareMathAlphabet{\mathpzc}{OT1}{pzc}{m}{it}

\newcommand{\R}{\ensuremath{\mathbb R}}

\newcommand{\T}{\ensuremath{\mathbb T}}
\newcommand{\Y}{\ensuremath{\mathbb Y}}
\newcommand{\Pp}{\ensuremath{\mathbb P}}
\newcommand{\1}{\mathbf{1}}

\newcounter{comptage}[part]
\newtheorem{lem}[comptage]{Lemma}
\newtheorem{theo}[comptage]{Theorem}
\newtheorem{defin}[comptage]{Definition}

\newtheorem{prop}[comptage]{Proposition}
\newtheorem{conj}[comptage]{Conjecture}

\newenvironment{theo_n}[2][Theorem]{\begin{trivlist}
\item[\hskip \labelsep {\bfseries #1}\hskip \labelsep {\bfseries #2}] \itshape}{\end{trivlist}}

\newenvironment{prop_n}[2][Proposition]{\begin{trivlist}
\item[\hskip \labelsep {\bfseries #1}\hskip \labelsep {\bfseries #2}] \itshape}{\end{trivlist}}

\newtheoremstyle{remarque}
  {3pt}
  {3pt}
  {}
  {}
  {\bf}
  {.}
  {.5em}
  {}
\theoremstyle{remarque}
\newtheorem{rem}[comptage]{Remark}

\author{\vspace{1cm} Antoine Lemenant \\
Université Paris XI \\
antoine.lemenant@math.u-psud.fr}

\title{On the homogeneity of global minimizers for the Mumford-Shah functional when $K$ is a smooth cone.}

\begin{document}
\maketitle

{\bf Abstract.} We show that if $(u,K)$ is a global minimizer for
the Mumford-Shah functional in $\R^N$, and if $K$ is a smooth
enough cone, then (modulo constants) $u$ is a homogenous function of
degree $\frac{1}{2}$. We deduce some applications in $\R^3$ as for instance
that an angular sector cannot be the singular set of a global
minimizer, that if $K$ is a  half-plane then $u$ is the corresponding cracktip function of two variables, or
that if $K$ is a cone that meets $S^2$ with an union of $C^\infty$ curvilinear
convex polygones, then it is a $\Pp$, $\Y$ or $\T$.



\tableofcontents

\section*{Introduction}

The functional of D. Mumford and J. Shah \cite{ms} was introduced to solve an image segmentation
problem. If $\Omega$ is an open subset of $\R^2$, for example a
rectangle, and $g\in L^\infty(\Omega)$ is an image, one can get a segmentation by minimizing
$$J(K,u):= \int_{\Omega \backslash K}|\nabla u|^2dx+\int_{\Omega\backslash K}(u-g)^2dx+H^{1}(K)$$
over all the admissible pairs $(u,K)\in
\mathcal{A}$ defined by
$$\mathcal{A}:=\{(u,K); \; K \subset \Omega \text { is closed } , \; u \in W^{1,2}_{loc}(\Omega \backslash K) \}.$$

 Any solution $(u,K)$ that
minimizes $J$ represents a ``smoother'' version of the image and
the set $K$ represents the edges of the image.

Existence of minimizers is  a well known result (see for instance
\cite{dcl}) using $SBV$ theory.

The question of regularity for the singular set $K$ of a minimizer
is more difficult. The following conjecture is currently still open.

\begin{conj}[Mumford-Shah]{\rm \cite{ms} } Let $(u,K)$  be a reduced minimizer for the functional $J$. Then
$K$ is the finite union of  $C^1$ arcs.
\end{conj}

The term ``reduced'' just means that we cannot find another pair $(\tilde u, \tilde K)$ such that $K \subset \tilde K$ and $\tilde u$ is
an extension of $u$ in $\Omega \backslash \tilde K$.

Some partial results are true for the conjecture. For instance it
is known that  $K$ is $C^1$ almost everywhere (see \cite{d1},
\cite{b} and \cite{afp1}). The closest to the conjecture is probably the result of A. Bonnet \cite{b}. He proves
that if $(u,K)$ is a minimizer, then every isolated connected component of $K$ is a finite union of $C^1$-arcs. The approach of A. Bonnet
is to use blow up limits. If $(u,K)$ is a minimizer in $\Omega$ and $y$ is a fixed point, consider the sequences  $(u_k,K_k)$ defined by
$$u_k(x)=\frac{1}{\sqrt{t_k}}u(y+t_kx), \quad K_k=\frac{1}{t_k}(K-y), \quad \Omega_k=\frac{1}{t_k}(\Omega-y).$$
When  $\{t_k\}$ tends to infinity, the sequence $(u_k,K_k)$ may tend to a pair $(u_\infty,K_\infty)$, and then  $(u_\infty,K_\infty)$ is
 called a Global Minimizer. Moreover,
A. Bonnet proves that if $K_\infty$ is connected, then $(u_\infty, K_\infty)$ is one of the list below :

$\bullet${\bf 1st case}: $ K_\infty=\varnothing$  and  $u_\infty$ is a constant.

$\bullet${\bf 2nd case}: $K_\infty$ is a line and $u_\infty$ is locally constant.

$\bullet${\bf 3rd case}: ``\emph{Propeller}'': $K_\infty$ is the union of  3
half-lines meeting with  $120$ degrees and $u_\infty$ is locally constant.

$\bullet${\bf 4th case}:
``\emph{Cracktip}'': $K_\infty=\{(x,0);x\leq 0\}$ and
$u_\infty(r\cos(\theta),r\sin(\theta))=
\pm\sqrt{\frac{2}{\pi}}r^{1/2}\sin\frac{\theta}{2}+C$, for $r>0$
and $|\theta|<\pi$ ($C$ is a constant), or a similar pair obtained by translation and rotation.

We don't know
whether the list is complete without the hypothesis that $K_\infty$ is connected. This would give a positive answer to the Mumford-Shah conjecture.

The Mumford-Shah functional was initially given in dimension $2$ but there is no restriction to define Minimizers for the analogous functional in
$\R^N$. Then we can also do some blow-up limits and try to think about what should be a global minimizer in $\R^N$. Almost nothing is known in this direction
and this paper can be seen as a very preliminary step. Let state some definitions.

\begin{defin} Let $\Omega \subset \R^N$,  $(u,K) \in  \mathcal{A}$ and $B$ be a ball such that $\bar B\subset \Omega$.
A competitor  for the pair  $(u,K)$ in the ball $B$ is a pair
$(v,L) \in \mathcal{A}$ such that
$$
\left.
\begin{array}{c}
u=v \\
K=L
\end{array}
\right\} \text{ in } \Omega \backslash  B
$$
and in addition if $x$ and $y$ are two points in
$\Omega \backslash ( B \cup K) $ that are separated by $K$ then
they are also separated by $L$.
\end{defin}

The expression ``be separated by $K$'' means that $x$ and $y$ lie
in different connected components  of $\Omega \backslash K$.

\begin{defin}\label{defms} A global minimizer in $\R^N$  is a pair $(u,K)\in \mathcal{A}$ (with $\Omega=\R^N$)
such that for every ball $B$  in $\R^N$ and every
competitor $(v,L)$ in $B$ we have
$$\int_{B \backslash K}|\nabla u|^2dx +H^{N-1}(K\cap B)\leq \int_{B\backslash L}|\nabla v|^2dx+H^{N-1}(L\cap B)$$
where $H^{N-1}$ denotes
the Hausdorff measure of dimension $N-1$.
\end{defin}

Proposition 9 on page 267 of \cite{d} ensures that any blow up limit of a minimizer for the Mumford-Shah functional in $\R^N$, is a
global minimizer in the sense of Definition \ref{defms}. As a beginning for the description of global minimizers in $\R^N$, we can firstly think
 about what should be a global minimizer in $\R^3$. If $u$ is locally constant, then $K$ is a minimal cone, that is, a set that
 locally minimizes the Hausdorff measure of dimension 2 in $\R^3$. Then by \cite{d3} we know that $K$ is a cone of type $\Pp$ (hyperplane), $\Y$
 (three half-planes meeting with 120 degrees angles) or of type $\T$ (cone over the edges of a regular tetraedron centered at the origin).
 Those cones became famous by the theorem of J. Taylor \cite{ta} which says that any minimal surface in $\R^3$ is locally $C^1$ equivalent to a
 cone of type $\Pp$, $\Y$ or $\T$.
\begin{center}
\includegraphics[width=4cm]{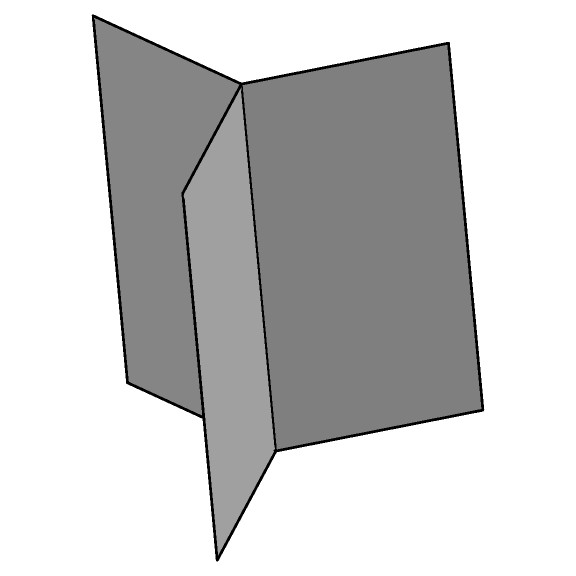} \hspace{2cm}
\includegraphics[width=4cm]{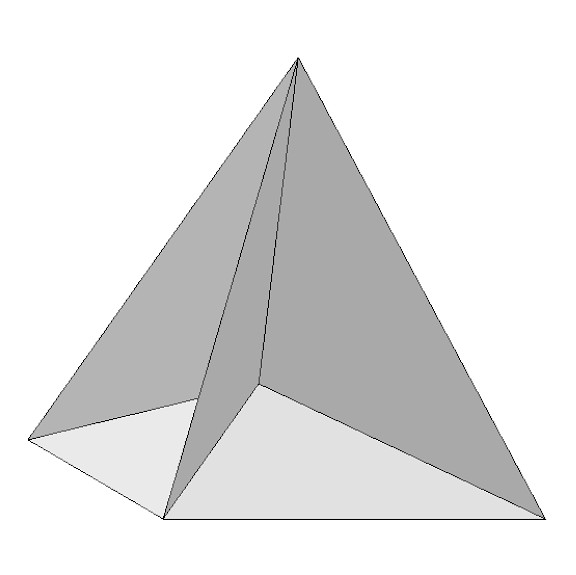}\\
Cones of type $\Y$ and $\T$ in $\R^3$.
\end{center}

To be clearer, this is a more precise  definition of $\Y$ and $\T$, as in  \cite{ddpt}.

\begin{defin}\label{prop} Define $Prop\subset \R^2$
by
$$Prop=\{(x_1,x_2);x_1 \geq 0, x_2=0\} $$
$$\hspace{4cm} \cup\{(x_1,x_2);x_1 \leq 0, x_2=-\sqrt{3}x_1\}$$
$$\hspace{6.5cm}\cup\{(x_1,x_2);x_1 \leq 0, x_2=\sqrt{3}x_1\}.$$
Then let $Y_0=Prop\times \R \subset \R^3.$ The spine of $Y_0$ is
the line $L_0=\{x_1=x_2=0\}$. A cone of type  $\Y$
is a set $Y=R(Y_0)$ where  $R$  is the composition of a
 translation and a rotation. The spine of $Y$ is then the line $R(L_0)$.
 \end{defin}

\begin{defin}\label{defT} Let $A_1=(1,0,0)$, $A_2=(-\frac{1}{3},\frac{2\sqrt{2}}{3},0)$,
$A_3=(-\frac{1}{3},-\frac{\sqrt{2}}{3},\frac{\sqrt{6}}{3})$, and
$A_4=(-\frac{1}{3},-\frac{\sqrt{2}}{3}, -\frac{\sqrt{6}}{3})$ the
four vertices of a regular tetrahedron centered at $0$. Let $T_0$
be the cone over the union of the $6$ edges $[A_i,A_j]$ $i\not
=j$. The spine of $T_0$  is the union of the four half lines
$[0,A_j[$. A cone of type $\T$ is a set
$T=R(T_0)$ where $R$ is the composition of a translation and a
rotation. The spine of $T$ is the image by $R$  of the spine of
$T_0$.
\end{defin}

So the pairs $(u,Z)$ where $u$ is locally constant and $Z$ is a minimal cone, are examples of global minimizers in
$\R^3$. Another global minimizer can be obtained with $K_\infty$ a half-plane, by setting $u:= Craktip\times \R$ (see \cite{d} section 76).
These examples are the only global minimizers in $\R^3$ that we know.

Note that if $(u,K)$ is a global minimizer in $\R^N$, then $u$ locally minimizes the Dirichlet integral
in $\R^N\backslash K$. As a consequence, $u$ is  harmonic in $\R^N \backslash K$. Moreover, if $B$ is a ball such that
 $K \cap B$ is regular enough, then the normal derivative of $u$ vanishes on $K \cap B$.

In this paper we wish to study global minimizers $(u,K)$ for which $K$ is a cone. It seems natural to think that any singular set of a
 global minimizer is a cone. But even if all known examples are cones, there is no proof of this fact. In addition, we will add some regularity on $K$.
 We denote by $S^{N-1}$ the unit sphere in $\R^N$ and, if $\Omega$ is a open set, $W^{1,2}(\Omega)$ is the Sobolev space. We will say that a domain $\Omega$ on $S^{N-1}$
 has a piecewise $C^2$ boundary, if the topological boundary of $\Omega$, defined by $\partial \Omega=\bar \Omega \backslash \Omega$, consists of
 an union of $N-2$ dimensional hypersurfaces of class $C^2$. This allows some cracks, i.e. when $\Omega$ lies in each sides of its boundary.
We will denote by $\tilde \Sigma$ the set of
 all the singular points of the boundary, that is
$$\tilde \Sigma:=\{x \in \partial \Omega; \forall r >0, B(x,r)\cap \partial \Omega \text{ is not a } C^2 \text{ hypersurface } \}.$$

\begin{defin} A smooth cone is a set $K$ of dimension $N-1$  in $\R^N$ such that $K$ is conical, centered at the origin,
and such that $S^{N-1}\backslash K$ is a domain with piecewise $C^2$ boundary.
Moreover we assume that the embedding $L^2(S^{N-1}\backslash K)\to W^{1,2}(S^{N-1}\backslash K)$ is compact.
Finally we suppose that we can strongly integrate by parts in $B(0,1)\backslash K$. More precisely,  denoting by $\Sigma$ the set of singularities
$$\Sigma:=\{tx ; (t,x) \in \R^+ \times \tilde \Sigma \},$$
we want that
 $$\int_{B(0,1)\backslash K}\langle \nabla u,\nabla \varphi \rangle=0$$
 for every harmonic function $u$ in $B(0,1)\backslash K$ with $\frac{\partial}{\partial n}u=0$ on $K \backslash \Sigma$, and for all $\varphi \in
 W^{1,2}(B(0,1)\backslash K)$ with vanishing trace on $S^{N-1}\backslash K$.
\end{defin}

\begin{rem} For instance, the cone over a finite union of $C^2$-arcs on $S^2$ is a smooth cone in $\R^3$. Another example in $\R^N$ is given by the union of admissible set of faces (as in Definition (22.2) of \cite{dau1}).
\end{rem}

Now this is the main result.

\begin{theo_n}{\ref{mainth}.}  Let $(u,K)$ be a global minimizer in $\R^N$. Assume that $K$ is a smooth cone.
Then there is a $\frac{1}{2}$-homogenous function $u_1$  such that $u-u_1$ is locally constant.
\end{theo_n}

As we shall see, this result implies that if $(u,K)$ is a global minimizer in $\R^N$, and if $K$ is a smooth cone other than a minimal cone,
then $\frac{3-2N}{4}$ is an eigenvalue for the spherical Laplacian in $S^{N-1}\backslash K$ with Neumann boundary conditions.
In section \ref{applications} we will give some applications about global minimizers in $\R^3$, using the estimates on the first eigenvalue
that can be found in \cite{da}, \cite{dau1} and \cite{kmr}. More precisely, we have :

\begin{prop_n}{\ref{app1}} Let $(u,K)$ be a global Mumford-Shah minimizer in
$\R^3$ such that $K$ is a smooth cone. Moreover, assume that $S^2\cap K $ is a union of convex curvilinear polygons with $C^\infty$
sides. Then  $u$ is locally constant and $K$ is a cone of type $\Pp$, $\Y$ or $\T$.
\end{prop_n}

Another consequence of the main result is the following.

\begin{prop_n}{\ref{cracktip}}  Let $(u,K)$ be a global Mumford-Shah minimizer in
$\R^3$ such that $K$ is a half plane.  Then $u$ is equal to a function of type $cracktip\times \mathbb{R}$, that is, in cylindrical coordinates,
$$u(r,\theta,z)=\pm \sqrt{\frac{2}{\pi}}r^{\frac{1}{2}}sin\frac{\theta}{2} +C$$
for $0<r< + \infty $, $-\pi< \theta < \pi$ where $C$ is a constant.
\end{prop_n}

Finally, we deduce two
 other consequences from Theorem \ref{mainth}. Let $(r,\theta,z)\in \R^+\times [-\pi,\pi]\times \R
$ be the cylindrical coordinates in $\R^3$. For all $\omega \in
[0,\pi]$ set
\begin{eqnarray}
\partial \Gamma_\omega:=\{(r,\theta,z) \in \R^3; \theta=-\omega \text{ or } \theta=\omega \}.\notag
\end{eqnarray}
and
\begin{eqnarray}
S_\omega:=\{(r,\theta,z) \in \R^3; z=0,\; r>0,\;\theta \in[-\omega, \omega]\;\}
\end{eqnarray}

Observe that $S_{0}$ is a half line, $S_{\frac{\pi}{2}}$, $\partial \Gamma_{0}$ and $\partial \Gamma_{\pi}$ are half-planes, and that
 $S_{\pi}$ and $\partial \Gamma_{\frac{\pi}{2}}$ are planes.

\begin{prop_n}{\ref{app3}} There is no global Mumford-Shah minimizer in
$\R^3$ such that $K$ is wing of type $\partial \Gamma_\omega$ with $\omega\not \in \{0, \frac{\pi}{2},\pi\}$.
\end{prop_n}

\begin{prop_n}{\ref{sect}}  There is no global Mumford-Shah minimizer in
$\R^3$ such that $K$ is an angular sector of type  $(u,S_\omega)$ for $\omega \not \in \{\frac{\pi}{2}, \pi\}$.
\end{prop_n}

{\bf Acknowledgements} : The author wishes to thank Guy David for having introduced him to the Mumford-Shah Functional, 
and for many helpful and interesting discussions on this subject.

\section{If $K$ is a cone then $u$ is  homogenous}

In this section we want to prove Theorem \ref{mainth}. Notice that this result is only
useful if the dimension $N\geq 3$. Indeed, in
dimension 2, if  $K$ is a cone then it is connected thus it is in
the list described in the introduction.

\subsection{Preliminary}

Let us recall a standard uniqueness result about energy minimizers.

\begin{prop} \label{stam1} Let $\Omega$ be an open and connected set of  $\R^N$ and let
$I\subset \partial \Omega$ be a hypersurface of class $C^\infty$. Suppose that $u$ and $v$ are two functions in $W^{1,2}(\Omega)$ such that
$u=v$ a.e. on $I$ (in terms of trace), solving
the minimizing problem
$$\min E(w):=\int_{\Omega}|\nabla w(x)|^2dx$$
 over all the functions $w\in W^{1,2}(\Omega)$ that are equal to $u$ and $v$ on $I$. Then
 $$u=v.$$
\end{prop}

{\bf Proof :} This comes from a simple convexity argument which can be found for instance in \cite{d}, but let us write the proof since
it is very short. By the parallelogram identity we have
\begin{eqnarray}
E(\frac{u+v}{2})=\frac{1}{2}E(u)+\frac{1}{2}E(v)-\frac{1}{4}E(u-v). \label{paral}
\end{eqnarray}
On the other hand, since $\frac{u+v}{2}$ is equal to $u$ and $v$ on $I$, and by minimality of $u$ and $v$ we have
$$E(\frac{u+v}{2})\geq E(u)=E(v).$$
Now by \eqref{paral} we deduce that $E(u-v)=0$ and since $\Omega$ is connexe, this implies that $u-v$ is a constant. But $u-v$ is equal to $0$ on $I$ thus
$u=v$. \qed

\begin{rem} The existence of a minimizer can also be proved using the convexity of $E(v)$.
\end{rem}

\subsection{Spectral decomposition}

The key ingredient to obtain the main result will be the spectral
theory of the Laplacian on the unit sphere. Since $u$ is harmonic, we will decompose $u$ as a sum of
homogeneous harmonic functions just like we usually use the classical spherical harmonics. The difficulty here comes
from the lack of regularity of $\R^N\backslash K$.

It will be convenient to work with connected sets. So let $\Omega$
be a connected component of $S^{N-1}\backslash K$, and let $A(r)$
be
$$A(r):=\{tx ; (x,t)\in  \Omega \times [0,r[ \; \}.$$
We also set
$$A(\infty):= \{tx ; (x,t)\in  \Omega \times \R^+\; \}.$$
 All the following results are using that the
embedding $W^{1,2}(\Omega)$ in $L^2(\Omega)$ is compact. Recall that this is the case by definition, since $K$ is a smooth cone.
Notice that for instance the cone property insures that the embedding is compact (see Theorem
6.2. p 144 of \cite{a}).

Consider the quadratic form
$$Q(u)=\int_{\Omega}|\nabla u(x)|^2dx$$
of domain $W^{1,2}(\Omega)$ dense into the Hilbert space
$L^{2}(\Omega)$. Since $Q$ is a positive and closed quadratic form
(see for instance Proposition 10.61 p.129 of \cite{lb}) there
exists a unique selfadjoint operator  denoted by $-\Delta_n$ of
domain $D(-\Delta_n)\subset W^{1,2}(\Omega)$ such that
$$\forall u \in D(-\Delta_n), \; \forall v \in W^{1,2}(\Omega),\quad
\int_{\Omega}\langle \nabla u, \nabla v\rangle=\int_{\Omega}
\langle -\Delta_nu,v \rangle.$$

\begin{prop}\label{defin laplacien neumann}
The operator
$-\Delta_n$ has a countably infinite discrete set of eigenvalues, whose eigenfunctions span  $L^{2}(\Omega)$.
\end{prop}

{\bf Proof :} The proof is the same as if $\Omega$ was a regular domain.  Consider the new quadratic form
$$\tilde{Q}(u):= Q(u)+\|u\|_2^2$$
with the same domain $W^{1,2}(\Omega)$. The form $\tilde{Q}$ has
the same properties than $Q$ and the associated operator is ${\rm
Id} - \Delta_n$. Moreover $\tilde{Q}$ is coercive. As a result,
the operator ${ \rm Id}-\Delta_n$ is  bijective and its inverse
goes from $L^{2}(\Omega)$ to $D(-\Delta_n)\subset
W^{1,2}(\Omega)$. By hypothesis the embedding of
$W^{1,2}(\Omega)$ into $L^2(\Omega)$ is compact. Thus the
resolvant $({ \rm Id}-\Delta_n)^{-1}$ is a compact operator,
and we conclude using the spectral theory of operators with a compact resolvant (see \cite{rs} Theorem XIII.64 p.245).  \qed

\begin{rem} The domain of $-\Delta_n$ is not known in general.
If $\Omega$ was smooth, then we could show that the domain is
exactly $D(-\Delta_n)=\{u\in W^{2,2}(\Omega); \frac{\partial
u}{\partial n}=0 \; \rm{on}\;
\partial \Omega\}$. Here, the boundary of $\Omega$ has some singularities so this result doesn't apply directly. But knowing exactly
the domain of $-\Delta_n$ will not be necessary for us.
\end{rem}

Now we want to study the link between the abstract operator $\Delta_n$ and the classical spherical
Laplacian $\Delta_S$ on the unit sphere. Recall that if we compute
the Laplacian in spherical coordinates, we obtain the following
equality
\begin{eqnarray}
\Delta = \frac{\partial^2 }{\partial
r}+\frac{N-1}{r}\frac{\partial}{\partial r}+\frac{1}{r^2}\Delta_S.
\label{forumlelaplacien}
\end{eqnarray}

\begin{prop} \label{regularite fonctions propres}
For every function  $f\in D(-\Delta_n)$ such that $-\Delta_n
f=\lambda f$ we have
\begin{eqnarray}
i) &\quad& f\in C^{\infty}(\Omega) \notag \\
ii)&\quad&  -\Delta_S f=-\Delta_n f=\lambda f \text { in }
\Omega \notag\\
iii) &\quad& \frac{\partial f}{\partial n} \text{ exists and is
equal to } 0 \text{ on } K \cap \overline{\Omega} \backslash
\Sigma \notag
\end{eqnarray}
\end{prop}
{\bf Proof :} Let $\varphi$ be a  $C^{\infty}$ function with
compact support in $\Omega$ and $f\in D(-\Delta_n)$. Then the
Green formula in the distributional sense gives
$$\int_{\Omega} \nabla f . \nabla \varphi  = \langle -\Delta_S f , \varphi \rangle$$
where the left and right brackets mean the duality in the
distributional sense. On the other hand, by definition of
$-\Delta_n$ and since $f$ is in the domain $D(-\Delta_n)$, we also
have
$$\int_{\Omega} \nabla f . \nabla \varphi  =\langle -\Delta_n
f ,\varphi \rangle$$ where this time the brackets mean the scalar
product in $L^2$. Therefore
$$ \Delta_n f =\Delta_S f \quad \text{in}\;
\mathcal{D}'(\Omega).$$ In other words, $-\Delta_S f = \lambda f$
in $\mathcal{D}'(\Omega)$. But now since $f\in W^{1,2}(\Omega)$,
by hypoellipticiy of the Laplacian we know that $f$ is $C^{\infty}$
and that  $-\Delta_S f=\lambda f$ in the classical sense. That
proves $i)$ and $ii)$. We even know by the elliptic theory that, since $K\backslash \Sigma$ is regular,
$f$ is regular at the boundary on $K\backslash \Sigma$.\\
 Now consider a ball  $B$  such that the intersection with $K\cap \overline{\Omega}$
 does not meet $\Sigma$.
 Assume that   $B$  is cut in two parts $B^+$ and $B^-$ by $K$, and that $B^+$ is one part in $\Omega$.
Possibly by modifying $B$ in a neighborhood of the intersection
with $K$, we can assume that the boundary of $B^+$ and $B^-$ is
$C^2$. The definition of $\Delta_n$ implies that for all function
$\varphi \in C^{2}(\bar \Omega)$ that vanishes out of $B^+$ we
have
$$\int_{B^+}\langle \nabla f , \nabla \varphi \rangle dx=
\int_{B^+}\langle-\Delta_n
 f, \varphi \rangle dx=\lambda\int_{B^+}\langle f , \varphi \rangle dx.$$
 On the other hand, integrating by parts,
\begin{eqnarray}
\int_{B^+}\langle \nabla f , \nabla \varphi \rangle dx&=&
\int_{B^+}\langle-\Delta_S f, \varphi \rangle+\int_{\partial
B^+}\frac{\partial u}{\partial n} \varphi \notag \\
&=&
 \lambda \int_{\partial B^+}\langle f, \varphi \rangle+\int_{\partial B^+}\frac{\partial f}{\partial n}
 \varphi\notag
\end{eqnarray}
thus
$$\int_{\partial B^+}\frac{\partial f}{\partial n} \varphi=0.$$
In other words the function $f$ is a weak solution of the mixed
boundary value problem
\begin{eqnarray}
-\Delta_S u=\lambda f & \text{ in }& B^{+} \notag \\
 u=f &\text{ on }& \partial B^{+}\backslash K \notag \\
 \frac{\partial u}{\partial n}=0 &\text{ on }& K\cap\partial
 B^{+}\notag
 \end{eqnarray}
Therefore, some  results from the elliptic theory imply
that $f$ is smooth in $B$ and is a strong solution (see \cite{taylorm}).\qed

Let us recapitulate what we have obtained. For all function $f\in
L^{2}(\Omega)$, there is a sequence of numbers $a_i$ such that
\begin{eqnarray}
f=\sum_{i=0}^{+\infty}a_i f_i \label{serie}
\end{eqnarray}
where the sum converges in  $L^{2}$. The functions $f_i$ are in
$C^{\infty}(\Omega)\cap W^{1,2}(\Omega)$, verify $-\Delta_S
f_i=\lambda_i f_i$ and $\frac{\partial f_i}{\partial n}=0$ on
$K\cap \overline{\Omega}\backslash \Sigma$. Moreover, we can normalize the $f_i$ in
order to obtain an orthonormal basis on $L^{2}(\Omega)$, in
particular we have the following Parseval formula \index{Parseval}
$$\|f\|_2^2=\sum_{i=0}^{+\infty}|a_i|^2. $$
Note that if $f$ belongs to the kernel of $-\Delta_n$ (i.e. is an
eigenfunction with eigenvalue  $0$), then
$$\langle \nabla f,\nabla f \rangle = \langle -\Delta_n f, f \rangle =0$$
and since $\Omega$ is connected that means that $f$ is a constant.
Thus $0$ is the first eigenvalue and the associated eigenspace has
dimension $1$. Then we can suppose that $\lambda_0=0$
and that all the $\lambda_i$ for $i > 0$ are positive.

We define the scalar product in $W^{1,2}(\Omega)$ by
$$\langle u,v \rangle_{W^{1,2}}:=\langle u,v \rangle_{L^2}+\langle \nabla u, \nabla v \rangle_{L^2}. $$

\begin{prop} \label{gradientsok} The family $\{f_i\}$ is orthogonal
in $W^{1,2}(\Omega)$. Moreover if $f\in W^{1,2}(\Omega)$ and if its
decomposition in $L^2(\Omega)$ is $f=\sum_{i=0}^{+\infty} a_i f_i$,
then the sum $\sum_{i=0}^{+\infty}|a_i|^2 \|\nabla f_i\|_{2}^{2}$
converges and
\begin{eqnarray}
\sum_{i=0}^{+\infty}|a_i|^2 \|\nabla f_i\|_{2}^{2}= \|\nabla
f\|_{2}^{2}. \label{parseval}
\end{eqnarray}
\end{prop}

{\bf Proof :} We know that $\{f_i\}$ is an orthogonal family in
$L^2(\Omega)$. In addition if $i\not = j$ then
\begin{eqnarray}
\int_{\Omega}\nabla f_i \nabla f_j&=&\int_{\Omega}-\Delta_n f_i
f_j \notag \\
&=& \lambda_i \int_{\Omega}f_i f_j \notag \\
&=& 0 \notag
\end{eqnarray}
thus $\{f_i\}$ is also orthogonal in $W^{1,2}(\Omega)$ and
$$\|f_i\|_{W^{1,2}}^2:=\|f_i\|_2^2+\|\nabla f_i\|_{2}^2=1+\lambda_i.$$
Consider now the orthogonal projection (for the scalar product of
$L^2$)
$$P_k:f \mapsto \sum_{i=0}^{k}a_i f_i.$$
The operator $P_k$ is the orthogonal projection on the closed
subspace  $A_k$ generated by $\{f_{0},...,f_k\}$. More precisely,
we are interested in the restriction of $P_k$ to the subspace
$W^{1,2}(\Omega)\subset L^2(\Omega)$. Also denote by
$\tilde{P}_k:W^{1,2}\to A_k$ the orthogonal projection on the same
subspace but for the scalar product of $W^{1,2}$. We want to show
that $P_k=\tilde{P}_k$. To prove this, it suffice to show that for all sets of coefficients $\{a_i\}_{i=1..k}$ and $\{b_i\}_{i=1..k}$,
$$\langle f-\sum_{i=0}^k a_i f_i , \sum_{i=0}^k b_i f_i\rangle_{W^{1,2}}=0.$$
Since we already have
$$\langle f-\sum_{i=0}^k a_i f_i , \sum_{i=0}^k b_i f_i\rangle_{L^2}=0,$$
all we have to show is that
$$\int_{\Omega}\langle \nabla f-\sum_{i=0}^k a_i \nabla f_i,\sum_{i=0}^k b_i \nabla f_i \rangle dx=0.$$
Now
\begin{eqnarray}
\int_{\Omega}\langle \nabla f-\sum_{i=0}^k a_i \nabla
f_i,\sum_{i=0}^k b_i \nabla f_i \rangle &=& \int_{\Omega}\langle
\nabla f,\sum_{i=0}^k b_i \nabla f_i \rangle - \sum_{i=0}^k a_ib_i \|\nabla f_i\|_2^2\notag \\
&=&\sum_{i=0}^k b_i \langle -\Delta_n f_i,f
\rangle_{L^2}-\sum_{i=0}^{k}a_ib_i \lambda_i \notag \\
&=&\sum_{i=0}^{k}a_ib_i \lambda_i \notag-\sum_{i=0}^{k}a_ib_i
\lambda_i \notag \\
&=&0 \notag
\end{eqnarray}
thus $P_k=\tilde{P}_k$ and therefore, by Pythagoras
$$\|P_{k}(f)\|_{W^{1,2}}^{2}\leq \|f\|_{W^{1,2}}^2.$$
By letting  $k$ tend to infinity we obtain
\begin{eqnarray}
\sum_{i=0}^{+\infty}a_i^2 \|\nabla f_i\|_{2}^{2}\leq \|\nabla
f\|_{2}^{2}. \label{ineq}
\end{eqnarray}
 From this inequality we deduce that the sum is absolutely converging
in $W^{1,2}(\Omega)$. Therefore, the sequence of partial sum $\sum_{i=0}^{K}a_i f_i$ is a Cauchy sequence for the norm $W^{1,2}(\Omega)$.  Thus, since
the sum  $\sum a_i  f_i $ already converges to $f$ in $L^2(\Omega)$, by uniqueness of the limit the sum converges to $f$ in $W^{1,2}(\Omega)$, so we deduce that
\eqref{ineq} is an equality  and the prove is over.\qed

Once we have a basis $\{f_i\}$ on $\Omega \subset S^{N-1}$, we consider for a certain $r_0>0$, the functions
$$h_i(x)=r_0^{\alpha_i}f_i\left(\frac{x}{r_0}\right)$$
defined on $r_0\Omega$. The exponent $\alpha_i$ is defined by
 \begin{eqnarray}
 \alpha_i=\frac{-(N-2)+\sqrt{(N-2)^2+4\lambda_i}}{2}. \label{defalphai}
 \end{eqnarray}
\begin{sloppypar}
The functions  $h_i$ form a basis of  $W^{1,2}(r_0\Omega)$.
Indeed, if $f\in W^{1,2}(r_0\Omega)$, then $f(r_0x)\in
W^{1,2}(\Omega)$ thus applying the decomposition on $\Omega$ we
obtain
\end{sloppypar}
$$f(r_0x)=\sum_{i=0}^{+\infty}b_i f_i(x)$$
thus
$$f(x)=\sum_{i=0}^{+\infty}a_i h_i(x)$$
with
\begin{eqnarray}
a_i=b_i r_0^{-\alpha_i}. \label{aibi}
\end{eqnarray} Notice that since
$\|h_i\|_{2}^2=r_0^{2\alpha_i+N-1}$ we also have
\begin{eqnarray}
\sum_{i=0}^{\infty}a_i^2
\|h_i\|^2_2=\sum_{i=0}^{\infty}a_i^2r_0^{2\alpha_i+N-1}
=\|f\|^2_{L^2(r_0\Omega)}<+\infty. \label{convhi}
\end{eqnarray}
Moreover, applying Proposition \ref{gradientsok} we have that
\begin{eqnarray}
\sum_{i=0}^{\infty}b_i^2 \|\nabla f_i\|_2^2 = \|\nabla f(r_0
x)\|_{2}^2 <+\infty. \label{gradhi}
\end{eqnarray}

We are now able to state our decomposition in $A(r_0)$.

\begin{prop}\label{decomp1} Let $K$ be a  smooth cone in $\R^N$, centered at the
origin and let $\Omega$ be a connected component of
$S^{N-1}\backslash K$. Then there exist some harmonic
homogeneous functions $g_i$, orthogonal in $W^{1,2}(A(1))$, such
that for every function $u\in W^{1,2}(A(1))$ harmonic in $A(1)$
with $\frac{\partial u}{\partial n}=0$  on $K\cap
A(1)\backslash \Sigma$, and for every $r_0\in ]0, 1[$, we
have that
$$u=\sum_{i=0}^{+\infty}a_i g_i \quad \text{ in } A(r_0)$$
where the $a_i$ do not depend on radius $r_0$ and are unique. The
sum converges in $W^{1,2}(A(r_0))$ and uniformly on all compact
sets of $A(1)$. Moreover
\begin{eqnarray}
\|u\|_{W^{1,2}(A(r_0))}^2=\sum_{i=0}^{+\infty}a_i^2\|g_i\|^2_{W^{1,2}(A(r_0))}.
\label{etoile}
\end{eqnarray}
\end{prop}

{\bf Proof  :}
 Since $u\in W^{1,2}(A(1))$ then for almost every
$r_0$ in $]0,1]$ we have that
$$u|_{r_0\Omega}\in W^{1,2}(r_0\Omega).$$
Thus we can apply the decomposition on $r_0\Omega$ and say that
$$u=\sum_{i=0}^{+\infty}a_i h_i \quad \text{ on }r_0\Omega.$$
Define $g_i$ by
$$g_i(x):=\|x\|^{\alpha_i}f_i\left(\frac{x}{\|x\|}\right)$$
where  $\alpha_i$ is defined by \eqref{defalphai}. Since the $f_i$ are eigenfunctions for $-\Delta_S$, we deduce from
 \eqref{forumlelaplacien} that
\begin{eqnarray}
\Delta g_i&=&\frac{\partial ^2}{\partial
r}g_i+\frac{N-1}{r}\frac{\partial}{\partial
r}g_i+\frac{1}{r^2}\Delta_Sg_i \notag \\
&=&\alpha_i(\alpha_i -1)r^{\alpha_i-2}f_i+\frac{N-1}{r}\alpha_i
r^{\alpha_i-1}f_i-r^{\alpha_i-2}\lambda_if_i \notag \\
&=&(\alpha_i^2+(N-2)\alpha_i-\lambda_i)r^{\alpha_i-2}f_i\notag \\
&=&0 \notag
\end{eqnarray}
by definition of $\alpha_i$, thus the $g_i$ are harmonic in $A(+\infty)$. Notice that the $g_i$ are orthogonal in $L^{2}(A(1))$ because they
are homogeneous and orthogonal in $L^{2}(\Omega)$. Note also that
$h_i$ is equal to $g_i$ on $r_0\Omega$. Moreover for all $0<r\leq
1$ we have
\begin{eqnarray}
\|g_i\|_{L^2(A(r))}^2&=& \int_{A(r)}|g_i|^2=\int_{0}^r\int_{\partial B(t)\cap A(1)}|g_i(w)|^2dwdt \notag \\
&=&\int_{0}^{r}\int_{\Omega}t^{N-1}|g_i(ty)|^2dydt =
\int_{0}^rt^{2\alpha_i+N-1}\int_{\Omega}|g_i(y)|^2dydt
\notag \\
&=&
\frac{r^{2\alpha_i+N}}{2\alpha_i+N}\|f_i\|^2_{L^2(\Omega)}=\frac{r^{2\alpha_i
+N}}{2\alpha_i+N}\leq 1. \label{normegi}
\end{eqnarray}

In the other hand, since the $f_i$ and their tangential gradients are
orthogonal in $L^{2}(\Omega)$, we deduce that the gradients of
$g_i$ are orthogonal in $A(1)$. Then, by a computation similar to
\eqref{normegi} we obtain for all $0<r\leq 1$
\begin{eqnarray}
\|\nabla g_i\|^2_{L^2(A(r))} &=& \int_{0}^r\int_{\partial B(t)\cap
A(1)}|\frac{\partial g_i}{\partial r}|^2
+|\nabla_\tau g_i|^2dwdt \notag \\
&=&\int_{0}^r\int_{\partial B(t)\cap A(1)}|\alpha_i
t^{\alpha_i-1}f_i(\frac{w}{t})|^2
+|t^{\alpha_i}\nabla_\tau f_i(\frac{w}{t})\frac{1}{t}|^2dwdt \notag \\
&=&\alpha_i^2\int_{0}^r t^{2(\alpha_i-1)}\int_{\partial B(t)\cap
A(1)}| f_i(\frac{w}{t})|^2dwdt+\int_{0}^r
t^{2(\alpha_i-1)}\int_{\partial B(t)\cap A(1)}
|\nabla_\tau f_i(\frac{w}{t})|^2dwdt \notag \\
&=&\alpha_i^2\int_{0}^r t^{2(\alpha_i-1)}\int_{\Omega}|
f_i(w)|^2t^{N-1}dw dt+\int_{0}^r t^{2(\alpha_i-1)}\int_{\Omega}
|\nabla_\tau f_i(w)|^2t^{N-1}dw dt \notag \\
&=&\alpha_i^2\frac{r^{2(\alpha_i-1)+N}}{2(\alpha_i-1)+N}\|f_i\|^2_{L^2(\Omega)} +
\frac{r^{2(\alpha_i-1)+N}}{2(\alpha_i-1)+N}\|\nabla_\tau f_i\|^2_{L^2(\Omega)} \notag \\
&=&\frac{r^{2(\alpha_i-1)+N}}{2(\alpha_i-1)+N}(\alpha_i^2+\lambda_i)\|f_i\|^2_{L^2(\Omega)}\notag \\
&\leq & Cr^{2\alpha_i}(\alpha_i^2+\lambda_i )
\label{estimgradients}
\end{eqnarray}
because $\|\nabla_\tau f_i\|_2^2=\lambda_i \|f_i\|_2^2$, $r\leq 1$ and $\alpha_i\geq 0$. Moreover the
constant $C$ depends on the dimension $N$ but does not depend on $i$.

We denote by $g$ the function defined in $A(\infty)$ by
$$g:=\sum_{i=0}^{+\infty}a_i g_i.$$
Then $g$ lies in  $L^{2}(A(r_0))$ because using
\eqref{normegi} and \eqref{convhi}
$$\|g\|_{L^2(A(r_0))}^2=\sum_{i=0}^{+\infty}|a_i|^{2}
\|g_i\|_{L^2(A(r_0))}^2\leq
\sum_{i=0}^{+\infty}|a_i|^2r_0^{2\alpha_i+N} <+\infty.$$ We want now
to show that  $g=u$.

$\bullet$ \emph{First step :} We claim that $g$ is harmonic in
$A(r_0)$. Indeed, since the $g_i$ are all harmonic  in
$A(r_0)$, the sequence of partial sums $s_k:=\sum_{i=0}^{k}a_i
g_i$ is a sequence of harmonic functions, uniformly bounded for
the $L^2$ norm in each compact set of $A(r_0)$. By the Harnack
inequality we deduce that the sequence of partial sums is
uniformly bounded for the uniform norm in each compact set. Thus
there is a subsequence that converges uniformly to a harmonic
function, which in fact is equal to $g$ by uniqueness of the
limit. Therefore,  $g$ is harmonic in $A(r_0)$.

$\bullet$ \emph{Second step :} We claim that $g$ belongs to
$W^{1,2}(A(r_0))$. Firstly, since $u\in W^{1,2}(r_0\Omega)$, by
\eqref{aibi} and \eqref{gradhi} we have that
\begin{eqnarray}
\sum_{i=0}^{+\infty}a_i^2r_0^{2\alpha_i} \| \nabla_\tau f_i
\|_{L^2(\partial B(0,1)\backslash K)}^2<+\infty
.\label{estimgradients2}
\end{eqnarray}
In addition,  since $\|\nabla_\tau f_i\|_2^2=\lambda_i\|f_i\|_{2}^2$
and $\|f_i\|_{2}=1$, we deduce
\begin{eqnarray}
\sum_{i=0}^{+\infty}a_i^2 r_0^{2\alpha_i}\lambda_i<+\infty
\label{estimgradients3}
\end{eqnarray}
and since $\alpha_i$ and $\lambda_i$ are linked by the formula
\eqref{defalphai} we also have that
\begin{eqnarray}
\sum_{i=0}^{+\infty}a_i^2 r_0^{2\alpha_i}\alpha_i^2<+\infty.
\label{estimgradients4}
\end{eqnarray}

Now, since $\sum a_i g_i$ converges absolutely on every compact set,
we can say that
$$\nabla g = \sum_{i=0}^{+\infty}a_i \nabla g_i$$
thus using \eqref{estimgradients}, \eqref{estimgradients3}, \eqref{estimgradients4}, and orthogonality,
\begin{eqnarray}
\|\nabla g\|^2_{L^{2}(A(r_0))} &=& \sum_{i=0}^{+\infty}
a_i^2 \|\nabla g_i\|_{L^2}^{2} \notag \\
& \leq& C\sum_{i=0}^{+\infty}a_i^2r_0^{2\alpha_i}(
\alpha_i^2+\lambda_i)<+\infty \notag.
\end{eqnarray}
Therefore, $g\in W^{1,2}(A(r_0))$. \\

$\bullet$ \emph{Third step :} We claim that $\frac{\partial
g}{\partial n}=0$ on  $K\cap \overline{A(r_0)} \backslash \Sigma$. We already
know that $\frac{\partial g_i}{\partial n}=0$ on $K\backslash
\Sigma$ (because the  $f_i$ have this property). We want to show that
$g$ is so regular that we can exchange  the order of $\frac{\partial}{\partial n}$ and $\sum$. So let  $x_0$ be
a point of  $K\cap \overline{A(r_0)}\backslash \Sigma$ and let $B$
be a neighborhood of $x_0$ in $\R^N$ that doesn't meet $\Sigma$ and
such that $K$ separates $B$ in two parts $B^+$ and $B^-$. Assume
that $B^+$ is a part in $A(r_0)$. The sequence of partial sums
$s_k:=\sum_{i=0}^k a_i g_i$ is a sequence of harmonic functions in
$B^+$. Since $\partial B^+\cap K$ is $C^2$ we can do a reflection
to extend $s_k$ in $B^-$. For all $k$, this new function $s_k$ is
the  solution of a certain elliptic equation whose operator become
from the composition of the Laplacian with the application that makes
$\partial B^+\cap K$ flat. Thus since $\sum a_i g_i$ converges absolutely
for the  $L^2$ norm, by the Harnack inequality  $\sum a_i g_i$  converges absolutely
for the uniform norm in a smaller neighborhood
$B'\subset B$ that still contains $x_0$. Thus $s_k$ converges to a
$C^1$ function denoted by $s$, which is equal to $g$ on
$B^+$. And since $\frac{\partial s_k}{\partial n}(x_0)=0$, by the
absolute convergence of the sum we can exchange the order of the derivative and the
symbol $\sum$  so we deduce that $\frac{\partial s}{\partial
n}(x_0)=0$. Finally, since $s$ is equal to  $g$ on $B^+$ we deduce
that $g$ is  $C^1$ at the boundary and $\frac{\partial
g}{\partial n}=0$ at $x_0$.

$\bullet$ \emph{Fourth step :} we claim that $g$ is equal to $u$
on $r_0\Omega$. Let $r$ be a radius such that $r<r_0$. Then the
function $x\mapsto g_r(x):= g(r\frac{x}{r_0})$ is well defined for
$x\in r_0\Omega$, and since the  $g_i$ are homogeneous we have
$$g(r\frac{x}{r_0})=\sum_{i=0}^{+\infty}a_ig_i(r\frac{x}{r_0})=
\sum_{i=0}^{+\infty}\left(\frac{r}{r_0}\right)^{\alpha_i}a_ig_i(x)=
\sum_{i=0}^{+\infty}\left(\frac{r}{r_0}\right)^{\alpha_i}a_ih_i(x).$$
We deduce that the function $x\mapsto g(\frac{r}{r_0}x)$ is in
$L^{2}(r_0\Omega)$ and its coefficients in the basis $\{h_i\}$ are
 $\{(\frac{r}{r_0})^{\alpha_i} a_i\}$. We want to show that
$\|g_r-u\|_{L^2(r_0\Omega)}$ tend to $0$. Indeed,  writing $u$ in
the basis $\{h_i\}$
$$u=\sum_{i=0}^{+\infty}a_i h_i,$$
we obtain
\begin{eqnarray}
\|g_r-u\|_2^2=
\sum_{i=0}^{+\infty}\left(\left(\frac{r}{r_0}\right)^{\alpha_i}-1\right)^2
a_i^2 \|h_i\|_2^2 \notag
\end{eqnarray}
which tends to zero when $r$ tends to $r_0$ by the dominated convergence theorem because $\left(\left(\frac{r}{r_0}\right)^{\alpha_i}-1\right)^2\leq 1$. Therefore,
 there is a subsequence for which $g_r$ tends to $u$
 almost everywhere. On the other hand, since $g$ is harmonic, the limit of $g_r$ exists and is equal to $g$.
That means that  $g$ tends to $u$ radially at almost every point
of  $r_0\Omega$.

$\bullet$ \emph{Fifth  step:} The functions $u$ and $g$ are
harmonic functions in $A(r_0)$, with finite energy, with a normal derivative
equal to zero on $K\cap \overline{A(r_0)}\backslash \Sigma$ and
that coïncide on $\partial A(r_0)\backslash K$. To show that
$u=g$ in $A(r_0)$ we shall prove that $g$ is an energy minimizer. Proposition  \ref{stam1} will then give
the uniqueness.

Let $\varphi \in W^{1,2}(A(r_0))\backslash K)$ have a vanishing trace on
$\partial B(0,r_0)$. Then, setting  $J(v):=
\int_{A(r_0)}|\nabla v|^2$ for  $v\in
W^{1,2}(A(r_0))$ we have
$$J(g+\varphi)=J(g)+\int_{A(r_0)}\nabla g \nabla \varphi + J(\varphi).$$
Now since  $g$ is harmonic with  Neumann condition on
$K\backslash \Sigma$ and since $\varphi$ vanishes on $r_0\Omega$, integrating by parts  we obtain
$$J(g+\varphi)=J(g)+J(\varphi).$$
Since $J$ is non negative and $g+\varphi$ describes all the
functions in $W^{1,2}(A(r_0))$ with trace equal to
$u$ on $r_0\Omega$, we deduce that $g$
minimizes $J$. We can do the same with $u$ thus  $u$ and $g$ are
two energy minimizers with same boundary conditions. Therefore, by
Proposition \ref{stam1} we know that $g=u$.

$\bullet$ \emph{Sixth  step :} The decomposition do not
depends on $r_0$. Indeed, let  $r_1$ be a second choice of radius.
Then we can do the same work as before to obtain a decomposition

$$u(x):=\sum_{i=0}^{+\infty}b_i g_i(x) \quad \text{ in }
B(0,r_1)\backslash K.$$

Now by uniqueness of the decomposition in $B(0,min(r_0,r_1))$ we
deduce that $b_i=a_i$ for all $i$.

In addition, $r_0$ was initially chosen almost everywhere in $]0,1[$. But
since the decomposition does not depend on the choice of radius, $r_0$ can
be chosen anywhere in  $]0,1[$, by choosing a radius almost everywhere in $]r_0,1[$.
 \qed

\begin{theo} \label{mainth} Let  $(u,K)$ be a global minimizer in $\R^N$ such that
 $K$ is a smooth cone. Then for each connected
component of $\R^N\backslash K$ there is a constant $u_k$ such
that $u-u_k$ is $\frac{1}{2}$-homogenous.
\end{theo}

{\bf Proof :} Let $\Omega$ be a connected component of
$\R^N\backslash K$. We apply the preceding proposition to $u$.
Thus
$$u(x)=\sum_{i=0}^{+\infty}a_i g_i(x) \quad \text{ in } A(r_0).$$ for a certain radius $r_0$ chosen in $]0,1[$. Let us prove that the same decomposition is true in $A(\infty)$. Applying Proposition \ref{decomp1} to the function $u_R(x)=u(R x)$
we know that there are some coefficients $a_i(R)$ such that
$$u_R(x)=\sum_{i=0}^{+\infty}a_i(R)g_i(x) \text{ in } A(r_0).$$ Now since $u_R(\frac{x}{R})=u(x)$ we can use the homogeneity
of the $g_i$ to identify the terms in $B(0,r_0)$ thus $a_i(R)=a_i
R^{\alpha_i}$. Now we fix $y=Rx$ and we obtain that
$$u(y)=\sum_{i=0}^{+\infty}a_i g_i(y) \text{ in } A(Rr_0).$$ Since $R$ is arbitrary the decomposition is true in
$A(\infty)$.

In addition for every radius  $R$ we know that
\begin{eqnarray}
\|\nabla u\|^2_{L^{2}(A(R))}=
\sum_{i=0}^{+\infty}a_i^2 \|\nabla g_i\|^2_{L^2(A(R))} \label{som}
\end{eqnarray}
 and since $g_i$ is $\alpha_i$-homogenous,
$$\|\nabla g_i\|_{L^2(A(R))}^2=R^{2(\alpha_i-1)+N}\|\nabla g_i\|_{L^2(A(1))}^2.$$
Now, since  $u$ is a global minimizer, a classical estimate on the gradient obtained by comparing $(u,K)$ with $(v,L)$ where
 $v=\1_{B(0,R)^c}u$ and $L=\partial B(0,R)\cup (K\backslash B(0,R))$ gives
that there is a constant  $C$ such that for all radius  $R$
$$\|\nabla u\|_{L^2(B(0,R)\backslash K)}^2\leq C R^{N-1}.$$
We deduce
 $$ \sum_{i=0}^{+\infty}a_i^2 R^{2(\alpha_i-1)+N}\|\nabla g_i\|^2_{L^2(A(1))}\leq CR^{N-1}.$$
Thus
 $$ \sum_{i=0}^{+\infty}a_i^2 R^{2\alpha_i-1}\|\nabla g_i\|^2_{L^2(A(1))}\leq C.$$
This last quantity is bounded when $R$ goes to infinity if and
only if $a_i=0$ whenever $\alpha_i> 1/2$. On the other hand, this quantity is bounded when $R$ goes to $0$,
if and only if $a_i=0$ whenever $0<\alpha_i< 1/2$. Therefore, $u-a_0$ is a finite sum of terms of degree
$\frac{1}{2}$.\qed

\begin{rem} In Chapter 65 of \cite{d}, we can find a
variational argument that leads to a formula in dimension $2$ that
links the radial and tangential derivatives of $u$. For all $\xi
\in K\cap
\partial B(0,r)$, we call $\theta_\xi \in [0,\frac{\pi}{2}]$
the non oriented angle between the tangent to $K$ at point  $\xi$
and the radius $[0,\xi]$. Then we have the following formula
$$\int_{\partial B(0,r)\backslash K}\left(\frac{\partial u}{\partial r}\right)^2dH^1=
\int_{\partial B(0,r)\backslash K}\left(\frac{\partial u}{\partial
\tau}\right)^2dH^1+\sum_{\xi \in K \cap \partial B(0,r)}\cos
\theta_\xi -\frac{1}{r}H^1(K\cap B(0,r)).$$ Notice that for a
global minimizer in $\R^2$ with $K$ a centered cone we find
\begin{eqnarray}
\int_{\partial B(0,r)\backslash K}\left(\frac{\partial u}{\partial
r}\right)^2dH^1= \int_{\partial B(0,r)\backslash
K}\left(\frac{\partial u}{\partial \tau}\right)^2dH^1.
\label{deriv}
\end{eqnarray}
Now suppose that $(u,K)$ is a global minimizer in $\R^N$ with $K$
a smooth cone centered at $0$. Then by Theorem \ref{mainth} we know
that $u$ is harmonic and $\frac{1}{2}$-homogenous. Its restriction
to the unit sphere is an eigenfunction for the spherical Laplacian
with Neumann boundary condition and associated to the eigenvalue
$\frac{2N-3}{4}$. We deduce that
$$\|\nabla_\tau u\|_{L^2(\partial B(0,1))}^2=\frac{2N-3}{4}\|u\|_{L^2(\partial
B(0,1))}^2.$$ On the other hand $$\frac{\partial u}{\partial
r}(x)=\frac{1}{2}\|x\|^{-\frac{1}{2}}u(\frac{x}{\|x\|})$$ thus
$$\|\frac{\partial u}{\partial r}\|_{L^2(\partial B(0,1))}^2=\frac{1}{4}\|u\|_{L^2(\partial B(0,1))}^2.$$
So
$$\|\nabla_\tau u\|_{L^2(\partial B(0,1))}^2=(2N-3)\|\frac{\partial u}{\partial r}\|_{L^2(\partial B(0,1))}^2.$$
In particular, for $N=2$ we have the same formula as \eqref{deriv}.
\end{rem}

\section{Some applications} \label{applications}

As it was claimed in the introduction, here is some few applications of Theorem \ref{mainth}.

\begin{prop} \label{app1} Let $(u,K)$ be a global minimizer in
$\R^3$ such that $K$ is a smooth cone. Moreover, assume that $S^2\cap K $ is a union of convex curvilinear polygons with $C^\infty$
sides. Then  $u$ is locally constant and $K$ is a cone of type $\Pp$, $\Y$ or $\T$.
\end{prop}

{\bf Proof :} In each polygon we know by Proposition 4.5. of \cite{da} that the smallest positive eigenvalue for the operator minus
Laplacian with Neumann boundary conditions is greater than or equal to $1$. Thus it cannot be $\frac{3}{4}$ and $u$ is locally constant.
Then $K$ is a minimal cone in $\R^3$ and we know from \cite{d3} that it is a cone of type $\Pp$, $\Y$ or $\T$. \qed

Let $(r,\theta,z)\in \R^+\times [-\pi,\pi]\times \R
$ be the cylindrical coordinates in $\R^3$. For every $\omega \in
[0,\pi]$ set
$$\Gamma_\omega:=\{(r,\theta,z) \in \R^3; -\omega<\theta<\omega \}$$
of boundary
$$\partial \Gamma_\omega:=\{(r,\theta,z)\in \R^3; \theta=-\omega \text{ or } \theta=\omega \}.$$
Consider $\Omega_\omega= \Gamma_\omega \cap S^2$ and let
$\lambda_1$ be the smallest positive eigenvalue of $-\Delta_S$
in $\Omega_\omega$ with Neumann conditions on $\partial
\Omega_\omega$. Then by Lemma 4.1. of \cite{da} we have that
$$\lambda_1=\min(2,\lambda_\omega)$$
where
$$\lambda_\omega=\left(\frac{\pi}{2\omega}+\frac{1}{2}\right)^2-\frac{1}{4}.$$

In particular for the cone of type  $\Y$, $\omega=\frac{\pi}{3}$ thus
$\lambda_1=2$.

Observe that for $\omega\not = \pi $, $\lambda_\omega\not = \frac{3}{4}$. So we get this following proposition.

\begin{prop}\label{app3} There is no global Mumford-Shah minimizer in
$\R^3$ such that $K$ is wing of type $\partial \Gamma_\omega$ with $\omega\not \in \{0,\frac{\pi}{2},\pi \}$.
\end{prop}

Another consequence of Theorem \ref{mainth} is the following. Let $P$ be the half plane
$$P:=\{(r,\theta, z) \in \R^3; \theta=\pi\}.$$

\begin{prop} \label{cracktip} Let $(u,K)$ be a global Mumford-Shah minimizer in
$\R^3$ such that $K=P$. Then $u$ is equal to $cracktip\times \mathbb{R}$, that is in cylindrical coordinates
$$u(r,\theta,z)=\pm \sqrt{\frac{2}{\pi}}r^{\frac{1}{2}}sin\frac{\theta}{2} +C$$
for $0<r< + \infty $ and $-\pi< \theta < \pi$.
\end{prop}

\begin{rem} In Section \ref{ann} we will give a second proof of Proposition \ref{cracktip}.
\end{rem}

\begin{rem} We already know that $u=cracktip\times \R$ is a global minimizer in $\R^3$ (see \cite{d}).
\end{rem}

To prove Proposition \ref{cracktip} we will use the following well known result.

\begin{prop}[\cite{dau1}, \cite{k}] \label{cracktip2} The smallest positive eigenvalue for
 $-\Delta_n$ in $S^{2}\backslash P$ is $\frac{3}{4}$, the corresponding eigenspace is of dimension 
1 generated by the restriction on $S^2$ of the following  function in cylindrical coordinates $$u(r,\theta,z)=r^{\frac{1}{2}}sin\frac{\theta}{2} $$
for $0<r< + \infty $ and $-\pi< \theta < \pi$.
\end{prop}

Now the proof of Proposition \ref{cracktip} can be easily deduce from Proposition \ref{cracktip2} and Theorem \ref{mainth}.

{\bf Proof of Proposition \ref{cracktip}:} If $(u,P)$ is a global minimizer,
we know that after removing a constant the restriction of $u$ to the unit sphere is an eigenfunction for $-\Delta_n$ in $S^{2}\backslash P$ associated to the eigenvalue
 $\frac{3}{4}$. Therefore, from Proposition \ref{cracktip2} we know that
$$u(r,\theta,z)=Cr^{\frac{1}{2}}sin\frac{\theta}{2}$$
so we just have to determinate the constant $C$. But by a well known argument about Mumford-Shah minimizers we prove that $C$ must be equal to
$\pm \sqrt{\frac{2}{\pi}}$ (see \cite{d} Section 61 for more details).\qed

Now set
$$S_\omega:=\{(r,\theta,0); r>0,\theta \in[-\omega,\omega]\}$$

\begin{prop}\label{sect} There is no global Mumford-Shah minimizer in
$\R^3$ such that $K$ is an angular sector of type  $(u,
S_\omega)$ for $0<\omega<\frac{\pi}{2}$ or
$\frac{\pi}{2}<\omega<\pi$.
\end{prop}
{\bf Proof :} According to Theorem \ref{mainth}, if $(u, S_\omega)$ is a global minimizer, then
 $u-u_0$ is a homogenous harmonic function of degree  $\frac{1}{2}$, thus its restriction to $S^2\backslash S_\omega$ is an eigenfunction for
$-\Delta_n$ associated to the eigenvalue $\frac{3}{4}$. Now if
$\lambda(\omega)$ denotes the smallest eigenvalue on $\partial B(0,1)\backslash S_\omega$, we know by Theorem 2.3.2. p.47 of \cite{kmr}
 that $\lambda(\omega)$ is non decreasing with respect to $\omega$. Since $\lambda(\frac{\pi}{2})=\frac{3}{4}$,
we deduce that for $\omega < \frac{\pi}{2}$, we have
\begin{eqnarray}
\lambda(\omega)\geq \frac{3}{4}. \label{inega}
\end{eqnarray}
 In  \cite{kmr} page 53 we can find the following asymptotic formula near $\omega=\frac{\pi}{2}$
\begin{eqnarray}
\lambda(\omega)=\frac{3}{4}+\frac{2}{\pi}\cos\omega+O(\cos^2\omega).
\label{dev}
\end{eqnarray}
this proves that  the case when \eqref{inega} is a equality only arises when $\omega=\frac{\pi}{2}$. Thus  such eigenfunction $u$
doesn't exist.

Consider now the case $\omega> \frac{\pi}{2}$. For $\omega=\pi$ there are tow connected components. Thus $0$ is an eigenvalue of multiplicity 2. The
second eigenvalue is equal to $2$. Therefore, for $\omega=\pi$ the spectrum is
$$0 \leq 0 \leq 2 \leq \lambda_3\leq ... \quad \quad \omega= \pi$$
By monotonicity, when $\omega$ decreases, the eigenvalues increase. Since the domain becomes connexe, $0$
become of multiplicity 1 thus the second eigenvalue become positive. The spectrum is now
$$0 \leq \lambda_1 \leq \lambda_2 \leq  ... \quad \quad \omega < \pi$$
with $\lambda_i \geq 2$ for $i \geq 2$. Thus the only eigenvalue that could be equal to $\frac{3}{4}$ is $\lambda_2$ which is increasing from
 from $0$ to $\frac{3}{4}$, reached for $\omega=\frac{\pi}{2}$. Now \eqref{dev} says that the increasing is strict near $\omega=\frac{\pi}{2}$.
Therefore there is no eigenvalue  equal to
$3/4$ for $\omega > \frac{\pi}{2}$ and there is no possible global minimizer. \qed


\section{Second proof of Propositions \ref{cracktip} and \ref{cracktip2}} \label{ann}

Here we want to give a second proof of Proposition \ref{cracktip}, without using Theorem \ref{mainth}, and which do not use
 Proposition \ref{cracktip2}. In a remark at the end of this section, we will briefly explain how to use this  proof of
Proposition \ref{cracktip} in order to obtain a new proof of Proposition \ref{cracktip2} as well.

Let assume that $K$ is a half plane in $\R^3$. We can suppose for instance that
\begin{eqnarray}
K=P:=\{x_2=0\}\cap\{x_1\leq 0\} \label{definp}
\end{eqnarray}
 We begin by studying the harmonic measure in
$\R^3\backslash P$.

Let  $B:=B(0,R)$ be a ball of radius $R$ and let $\gamma$ be the trace operator
on $\partial B(0,R) \backslash P$. We denote by $T$ the image of $W^{1,2}(B\backslash K)$ by $\gamma$. We also denote by
$C_{b}^{0}(\partial B \backslash K)$ the set of continuous and bounded functions on $\partial B(0,1)\backslash P$. Finally set
$A:=T\cap C_{b}^{0} $. Obviously $A$ is not empty. To every function $f\in A $, Proposition 15.6. of \cite{d} associates a
unique energy minimizing function $u\in
W^{1,2}(B\backslash K)$ such that $\gamma(u)=f$ on $\partial B\backslash P$. Since $u$ is harmonic we know that it is $C^\infty$ in $B\backslash K$.
Let $y\in B\backslash K$ be a fixed point and consider the linear form $\mu_y$ defined by
\begin{eqnarray}
 \mu_y: A &\to&\R \label{defmuy} \\
 f&\mapsto&u(y). \notag
\end{eqnarray}

By the maximum principle for energy minimizers, we know that for all  $f\in A$
we have
$$| \mu_y(f)|\leq \|f\|_{\infty}$$
thus $\mu_y$ is a continuous linear form on
$A$ for the norm $\|\;\|_\infty$. We identify $\mu_y$ with its representant
in the dual space of $A$ and we call it  \emph{harmonic measure}.

Moreover, the harmonic measure is positive. That is, if $f\in
A$ is a non negative function, then  (by the maximum principle)
$\mu_y(f)$ is non negative. By positivity of
$\mu_y$, if $f \in A$ is a non negative function and $g \in A$
is such that  $fg\in A$, then since $(\|g\|_\infty+g)f$ and $(\|g\|_\infty-g)f$ are two non negative functions of
 $A$ we deduce that
\begin{eqnarray}
|\langle fg , \mu_y \rangle|\leq \|g\|_{\infty}\langle f,\mu_y
\rangle . \label{majoration}
\end{eqnarray}

Now here is an estimate on the measure  $\mu^R_y$.

\begin{lem}\label{rylem1} There is a dimensional constant $C_N$ such that the following holds. Let $R$ be a positive radius. For $0<\lambda<\frac{R}{2}$ consider the spherical domain
\begin{eqnarray}
\mathcal{C}_\lambda:=\{x \in \mathbb{R}^{3}\; ; \;|x|=R \textrm{ and } d(x,P)\leq \lambda\}. \notag
\end{eqnarray}
Let $\varphi_{\lambda} \in C^\infty(\partial B(0,R))$ be a function between $0$ and $1$, that
 is equal to $1$ on
$\mathcal{C}_\lambda$ and $0$ on $\partial B(0,R)\backslash \mathcal{C}_{2\lambda}$ and that
 is symmetrical with respect to $P$. Then for every
 $y \in B(0,\frac{R}{2})\backslash P$ we have
$$\mu^R_y(\varphi_{\lambda})\leq C\frac{\lambda}{R}.$$
\end{lem}
{\bf Proof :} Since $\varphi_{\lambda}$ is continuous and symmetrical with respect to  $P$, by the reflection principle, its harmonic
extension  $\varphi$ in $B(0,R)$ has a normal derivative equal to zero on $P$ in the interior of $B(0,R)$. Moreover
$\varphi_\lambda$ is clearly in the space $A$. Thus by definition of $\mu_y$,
$$\varphi(y)=\langle \varphi_{\lambda} ,  \mu^R_y \rangle .$$
On the other hand, since $ \varphi_{\lambda}$ is continuous on the entire sphere, we also have the formula with the classical Poisson kernel
\begin{eqnarray}
\varphi(y)=\frac{R^2-|y|^2}{N\omega_N R}\int_{\partial B_R}\frac{
\varphi_\lambda(x)}{|x-y|^{3}}ds(x) \notag
\end{eqnarray}
with $\omega_N$ equal to the measure of the unit sphere.
In other words
$$\mu^R_y(\varphi_\lambda)=\frac{R^2-|y|^2}{N\omega_N R}\int_{\partial B_R}\frac{
\varphi_\lambda(x)}{|x-y|^{3}}ds(x).$$ For $x\in \partial B_R$ we have
$$\frac{1}{2}R\leq |x|-|y|\leq |x-y|\leq |x|+|y|\leq \frac{3}{2}R.$$
We deduce that
$$\mu^R_y(\varphi_\lambda)\leq C_N \frac{1}{R^2}\int_{\mathcal{C}_{2\lambda}}ds.$$
Now integrating by parts,
\begin{eqnarray}
\int_{\mathcal{C}_\lambda}ds&=&2\int_{0}^\lambda 2\pi \sqrt{R^2-w^2}dw \notag \\
&=&4\pi\frac{\lambda}{2}\sqrt{R^2-\lambda^2}+R^2 \arcsin (\frac{\lambda}{R}) \notag\\
&\leq& C R\lambda \notag
\end{eqnarray}
because $\arcsin(x)\leq \frac{\pi}{2} x$. The proposition follows. \qed

Now we can prove the uniqueness of $cracktip \times \R$.

{\bf Second Proof of Proposition \ref{cracktip} :} Let us show that $u$ is vertically constant. Let $t$ be a positive real. For $x=(x_1,x_2,x_3)\in
\R^3$ set $x_t:=(x_1,x_2,x_3+t)$. We also set
$$u_t(x):=u(x)-u(x_t).$$
Since $u$ is a function associated to a global minimizer, and since $K$ is regular, we know that for all $R>0$, the restriction of  $u$
to the sphere $\partial B(0,R)\backslash K$ is continuous and bounded on
$\partial B(0,R)\backslash K$ with finite limits
on each sides of $K$. It is the same for
$u_t$. Thus for all $x\in \R^3\backslash P$ and for all
$R>2\|x\|$ we can write
$$u_t(x):= \langle u_t|_{\partial B(0,R)\backslash P}, \mu^R_x \rangle$$
where $\mu_x$ is the harmonic measure defined in \eqref{defmuy}. We want to prove that for $x \in \R^3 \backslash P$,
$\langle u_t|_{\partial B(0,R)\backslash P}, \mu^R_x \rangle$ tends to $0$ when  $R$ goes to infinity.
This will prove that $u_t=0$.

So let  $x\in \R^3\backslash P$ be fixed. We can suppose that $R>100(\|x\|+t)$.  Let $\mathcal{C}_\lambda$ and
 $\varphi_\lambda$ be as in Lemma \ref{rylem1}. Then write
\begin{eqnarray}
u_t(x)=  \langle u_t|_{\partial B(0,R)\backslash
P}\varphi_\lambda, \mu^R_x \rangle + \langle u_t|_{\partial
B(0,R)\backslash P}(1-\varphi_\lambda), \mu^R_x \rangle .\notag
\end{eqnarray}
Now by a standard estimate on Mumford-Shah minimizers (that comes from Campanato's Theorem, see \cite{afp} p. 371) we have for all $x \in \R^N
\backslash P$,
$$|u_t(x)|\leq C\sqrt{t}.$$
Then, using Lemma  \ref{rylem1} we obtain
\begin{eqnarray}
| \langle u_t|_{\partial B(0,R)\backslash P}\;\varphi_\lambda\;,\;
\mu^R_x \rangle| \leq C\sqrt{t}\frac{\lambda }{R}  \notag.
\end{eqnarray}
On the other hand, for the points $y$ such that
$d(y,P)\geq \lambda$, since $\tilde u: u(.)-u(y)$ is harmonic
in $B(y,d(y,P))$  we have, by a classical estimation on harmonic functions (see the introduction of \cite{gt})
$$|\nabla \tilde u(y)|\leq C\frac{1}{d(y,P)}\|\tilde u\|_{L^\infty(\partial B(y,\frac{1}{2}d(y,P)))}.$$
Now using Campanato's Theorem again we know that
$$\|\tilde u\|_{L^\infty(\partial B(y,\frac{1}{2}d(y,P)))}\leq
Cd(y,P)^{\frac{1}{2}}$$ thus
$$|\nabla  u(y)|\leq
C\frac{1}{d(y,P)^{\frac{1}{2}}}$$ and finally by the mean value theorem we deduce that for all
the points $y$ such that
$d(y,P)\geq \lambda$,
$$|u_t(y)| \leq C\sup_{z \in [y,y_t]}|\nabla u(z)|.|y-y_t|\leq t\frac{1}{\lambda^\frac{1}{2}}.$$
Therefore,
\begin{eqnarray}
| \langle u_t|_{\partial B(0,R)\backslash P}(1-\varphi_\lambda),
\mu^R_x \rangle| \leq Ct\frac{1 }{\lambda^{\frac{1}{2}}}. \notag
\end{eqnarray}
So
$$|u_t(x)|\leq C\sqrt{t}\frac{\lambda}{R}+Ct\frac{1 }{\lambda^{\frac{1}{2}}}$$
thus by setting $\lambda= R^{\frac{1}{2}}$ and by letting
$R$ go to $+\infty$ we deduce that $u_t(x)=0$ thus $z\mapsto
u(x,y,z)$ is constant.

Now we fix $z_0=0$ and we introduce $P_0:=P\cap \{z=0\}$.
We want to show that $(u(x,y,0),P_0)$ is a global minimizer in
$\mathbb{R}^2$. Let $(v(x,y),\Gamma)$ be a competitor for $u(x,y,0)$  in the 2-dimensional ball   $B$ of radius
$\rho$. Let $\mathcal{C}$ be the cylinder
$\mathcal{C}:=B\times[-R,R]$. Define $\tilde v$ and $\tilde
\Gamma$ in $\mathbb{R}^3$ by
$$\tilde v(x,y,z)= \left\{
\begin{array}{cc}
v(x,y) &\text{ if } (x,y,z)\in \mathcal{C} \\
u(x,y,z) &\text{ if } (x,y,z) \not \in \mathcal{C}
\end{array}\right.
$$

$$
\tilde \Gamma := (\mathcal{C}\cap[\Gamma \times [-R,R]])\cup
(P\backslash \mathcal{C})\cup (B\times \{\pm R\}).
$$
It is a topological competitor because
$\R^3\backslash P$ is connected (thus $P$ doesn't separate any points). Now finally let $\tilde B$ be a ball that contains
$\mathcal{C}$. Then $(\tilde v, \tilde \Gamma)$ is a competitor for
 $(u,P)$ in $\tilde B$. By minimality we have :
$$\int_{\tilde B}|\nabla u|^{2}+H^{2}(P \cap \tilde B) \leq \int_{\tilde B}|\nabla \tilde v|^{2}
+H^{2}(\tilde \Gamma \cap \tilde B).$$ In the other hand $u$ is equal to  $\tilde
v$ in $\tilde B \backslash \mathcal{C}$ and it is the same for $\Gamma$ and
$\tilde \Gamma$. We deduce
$$\int_{\mathcal{C}}|\nabla u|^{2}dxdydz+H^{2}(P \cap \mathcal{C}) \leq \int_{\mathcal{C}}|\nabla \tilde v|^{2}
dxdydz+H^{2}(\tilde \Gamma \cap \mathcal{C}).$$ Now, since
 $u$ and $\tilde v$ are vertically constant, $\nabla_z u =
\nabla_z \tilde v =0$, and $\nabla_x u$, $\nabla_y u$ are also constant with respect to the variable $z$ (as for $\tilde v$). Thus
$$2R\int_{B}|\nabla u(x,y,0)|^{2}dxdy+H^{2}(P \cap \mathcal{C}) \leq 2R\int_{B}|\nabla v(x,y)|^{2}dxdy +H^{2}(\tilde \Gamma \cap \mathcal{C}).$$
To conclude we will use the following lemma.

\begin{lem} \label{rotyi} If $\Gamma$ is rectifiable  and contained in a plane $Q$ then
$$H^{2}(\Gamma \times [-R,R])=2R H^{1}(\Gamma).$$
\end{lem}
{\bf Proof :} We will use the coarea formula (see Theorem
2.93 of \cite{afp}). We take $f:\mathbb{R}^{3}\to \mathbb{R}$ the orthogonal projection on the coordinate orthogonal to
 $Q$. By this way, if $E:=\Gamma \times [-R,R]$, we have $E\cap
f^{-1}(t)=\Gamma $ for all $t\in [-R,R]$. $E$ is rectifiable
(because $\Gamma$ is by hypothesis). So we can apply the coarea formula. To do this we have to calculate the jacobian $c_k d^{E}f_x$.
By construction, the approximate tangente plane in each point of $E$
is orthogonal to $Q$. We deduce that if $T_x$ is a tangent plane, then there is a basis of $T_x$
$(\overrightarrow{b_1},\overrightarrow{b_2})$ such that
$\overrightarrow{b_1}$ is orthogonal to $Q$. Since the function $f$
is the projection on $\overrightarrow{b_1}$, and its derivative as well
 (because $f$ is linear ) we obtain that the matrix of $d^{E}f_x:T_x \to
\mathbb{R}$ in the basis
$(\overrightarrow{b_1},\overrightarrow{b_2})$ is
$$d^{E}f_x=(1,0)$$
thus
$$c_k d^{E}f_x=\sqrt{det[(1,0).^t(1,0)]}=1.$$
Therefore
$$H^{2}(E)=\int_{-R}^{R}H^{1}(\Gamma)=2RH^{1}(\Gamma).\qed$$

Here we can suppose that $\Gamma$ is rectifiable. Indeed, the definition of Mumford-Shah minimizers is equivalent if we only 
allow rectifiables competitors. This is because the jump set of a $SBV$ function is rectifiable and in  \cite{dcl} it is 
proved that the relaxed functional on
the $SBV$ space has same minimizers.

So we have
$$2R\int_{B}|\nabla u(x,y,0)|^{2}dxdy+2R H^{1}(P \cap B) \leq 2R\int_{B}|\nabla v(x,y)|^{2}dxdy +
2RH^{1}(\Gamma \cap B)+H^2(B\times\{\pm R\}).$$ Then, dividing by $2R$,
$$\int_{B}|\nabla u(x,y,0)|^{2}dxdy+ H^{1}(P \cap B) \leq \int_{B}|\nabla v(x,y)|^{2}dxdy +
H^{1}(\Gamma \cap B)+\pi \frac{\rho^2}{R}$$ thus, letting
$R$ go to infinity,
$$\int_{B}|\nabla u(x,y,0)|^{2}dxdy+ H^{1}(P \cap B) \leq \int_{B}|\nabla v(x,y)|^{2}dxdy +
H^{1}(\Gamma \cap B).$$
This last inequality proves that $(u(x,y,0),P_0)$ is a global minimizer in $\mathbb{R}^{2}$, and since $P_0$ is a half-line, $u$ is a
$cracktip$.\qed

\begin{rem} Using a similar argument as the preceding proof, we can show that the first eigenvalue
for $-\Delta$ in $S^2\backslash P$ with Neumann boundary conditions (where $P$ is still a half-plane), is equal
to $\frac{3}{4}$. Moreover we can prove that the eigenspace is of dimension $1$,  generated by a function of type $cracktip\times \R$,
thus we have a new proof of Proposition \ref{cracktip2}. The argument is to take an eigenfunction $f$ in $S^2 \backslash P$, then to consider
 $u(x):=\|x\|^{\alpha}f(\frac{x}{\|x\|})$ with a good coefficient $\alpha\in ]0,\frac{1}{2}]$ that makes $u$ harmonic.  Finally we use the same sort of estimates on the harmonic measure to prove that $u$ is vertically constant. Thus we have reduced the problem in  dimension 2 and we conclude using that we know the eigenfunctions on the circle. A detailed proof is done in \cite{l1}.
\end{rem}

\section{Open questions}

As it is said in the introduction, this paper is a very short step in the discovering of all the global minimizers in $\R^N$.
 This final goal seems rather far but nevertheless some open questions might be accessible in a more reasonable time. All the following questions were 
pointed out by Guy David in  \cite{d}, and unfortunately they are still open after this paper.

$\bullet$ Suppose that $(u,K)$ is a global minimizer in $\R^N$. Is it true that $K$ is conical ?

$\bullet$ Suppose that $(u,K)$ is a global minimizer in $\R^N$, and $K$ is a cone. Is it true that $\frac{3-2N}{4}$
 is the smallest eigenvalue of the Laplacian on $S^{N-1}\backslash K$ ?

$\bullet$ Suppose that $(u,K)$ is a global minimizer in $\R^3$, and suppose that $K$ is contained in a plan (and not empty). 
Is it true that $K$ is a plane or a half-plane ?

$\bullet$ Could one found an extra global minimizer in $\R^3$ by blowing up the minimizer described in section 76.c. of \cite{d} (see also \cite{mer})?

One can find other open questions on global minimizers in the last page of \cite{d}.


\bibliographystyle{plain}
\bibliography{biblio}

ADDRESS :

Antoine LEMENANT \\
e-mail : antoine.lemenant@math.u-psud.fr

Université Paris XI\\
Bureau 15 Bâtiment 430 \\
ORSAY 91400 FRANCE

Tél: 00 33 169157951

\end{document}